\documentclass[12pt]{article}

\usepackage[T1]{fontenc}

\usepackage{amsmath, graphicx, color, amsthm, centernot, amssymb}
\usepackage[font=small,labelfont=bf,width=12.5cm]{caption}
\usepackage{algorithmic, algorithm}
\usepackage{subcaption}

\usepackage{fullpage}

\newcommand{\vc}[1]{\ensuremath{\vcenter{\hbox{#1}}}}

\newtheorem{proposition}{Proposition}
\newtheorem{conjecture}{Conjecture}

\newtheorem{problem}{Problem}

\if10     
\usepackage[mathlines]{lineno}
\newcommand*\patchAmsMathEnvironmentForLineno[1]{%
  \expandafter\let\csname old#1\expandafter\endcsname\csname #1\endcsname
  \expandafter\let\csname oldend#1\expandafter\endcsname\csname end#1\endcsname
  \renewenvironment{#1}%
     {\linenomath\csname old#1\endcsname}%
     {\csname oldend#1\endcsname\endlinenomath}}%
\newcommand*\patchBothAmsMathEnvironmentsForLineno[1]{%
  \patchAmsMathEnvironmentForLineno{#1}%
  \patchAmsMathEnvironmentForLineno{#1*}}%
\AtBeginDocument{%
\patchBothAmsMathEnvironmentsForLineno{equation}%
\patchBothAmsMathEnvironmentsForLineno{align}%
\patchBothAmsMathEnvironmentsForLineno{flalign}%
\patchBothAmsMathEnvironmentsForLineno{alignat}%
\patchBothAmsMathEnvironmentsForLineno{gather}%
\patchBothAmsMathEnvironmentsForLineno{multline}%
}
\linenumbers
\fi

\newcommand{\chiC}{\chi_{\mathrm{fum}}}

    {\noindent \emph{Proof.} {}{#1}{}}{$~$\hfill $~\blacklozenge$ \vspace{0.2cm}}

\title{A counterexample to a conjecture on facial unique-maximal colorings}
\author
{
	Bernard Lidick\'y\thanks{Department of Mathematics, Iowa State University, USA.
		E-Mail: \texttt{lidicky@iastate.edu}  } \and		
	Kacy Messerschmidt\thanks{Department of Mathematics, Iowa State University, USA.
		E-Mail: \texttt{kacymess@iastate.edu}} \and	
  	Riste \v{S}krekovski\thanks{Faculty of Information Studies, Novo mesto
  		\& University of Ljubljana, Faculty of Mathematics and Physics
  		\& University of Primorska, FAMNIT, Koper, Slovenia.
   		E-Mail: \texttt{skrekovski@gmail.com}}		
}

\begin{document}
\maketitle

{
	\abstract
	{
		A facial unique-maximum coloring of a plane graph is a proper vertex coloring by natural numbers where
		on each face $\alpha$ the maximal color appears exactly once on the vertices of $\alpha$.
		Fabrici and G\"{o}ring~\cite{FabGor16} proved that six colors are enough for any plane graph
		and conjectured that four colors suffice. This conjecture is a strengthening of the Four Color theorem.
		Wendland~\cite{Wen16} later decreased the upper bound from six to five.
		In this note, we disprove the conjecture by giving an infinite family of counterexamples.
		Thus we conclude that facial unique-maximum chromatic number of the sphere is five.		
	}

}

\section{Introduction}

We call a graph \emph{planar} if it can be embedded in the plane without crossing edges
and we call it \emph{plane} if it is already embedded in this way.
A \emph{coloring} of a graph is an assignment of colors to vertices.
A coloring is \emph{proper} if adjacent vertices receive distinct colors.
A proper coloring of a graph embedded on some surface,
where colors are natural numbers and every face has a unique vertex colored with a maximal color,
is called a \emph{facial unique-maximum coloring}, or \emph{FUM-coloring} for short.
The minimum $k$ such that a graph $G$ has a FUM-coloring using the colors $\{ 1, 2, \ldots, k \}$ is called the
\emph{facial unique-maximum chromatic number} of $G$ and is denoted $\chiC (G)$.

The cornerstone of graph colorings is the Four Color Theorem stating that every planar graph
can be properly colored using at most four colors~\cite{AppHak76}.
Fabrici and G\"{o}ring~\cite{FabGor16} proposed the following strengthening of the Four Color Theorem.

\begin{conjecture}[Fabrici and G\"{o}ring]\label{conj}
	\label{conj:plane4}
	If $G$ is a plane graph, then $\chiC (G) \leq 4$.
\end{conjecture}

When stating the conjecture, Fabrici and G\"{o}ring~\cite{FabGor16} proved that  $\chiC (G) \leq 6$ for every plane graph $G$.
Promptly, this coloring was considered by others.
Wendland~\cite{Wen16} decreased the upper bound to $5$ for all plane graphs.
Andova, Lidick\'y, Lu\v{z}ar, and \v{S}krekovski~\cite{1708.00094} showed that $4$ colors suffice for outerplanar graphs and for subcubic plane graphs.
Wendland~\cite{Wen16} also considered the list coloring version of the problem, where he was able
to prove the upper bound 7 and conjectured that lists of size 5 are sufficient.
Edge version of the problem was considered by Fabrici, Jendrol', and Vrbjarov\'a~\cite{FabJenVrb15}.
For more results on facially constrained colorings, see a recent survey written by Czap and Jendrol' \cite{CzaJen16}.

In this note we disprove Conjecture~\ref{conj}.

\begin{proposition}\label{prop}
	\label{prop:plane4}
	There exists a plane graph $G$ with $\chiC(G) > 4$.
\end{proposition}

\begin{figure}
	\centering
	\includegraphics{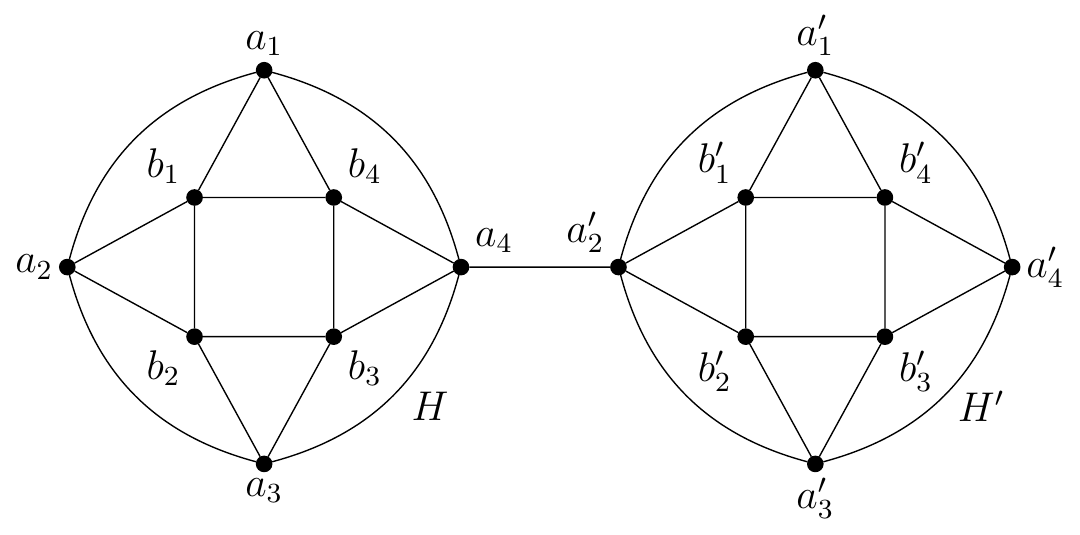}
	\caption{A counterexample to Conjecture~\ref{conj:plane4}.}
	\label{fig:plane41}
\end{figure}

\begin{proof}
	Let $G$ be the graph depicted in Figure~\ref{fig:plane41}.
	It consists of the induced graph $H$ on the vertex set $\{ a_1, a_2, a_3, a_4, b_1, b_2, b_3, b_4 \}$,
	$H'$ (an isomorphic copy of $H$), and the edge $a_4a_2'$ connecting them.
	Suppose for contradiction that $G$ has a FUM-coloring with the colors in $\{ 1, 2, 3, 4 \}$.
	The color $4$ is assigned to at most one vertex in the outer face of $G$,
	so by symmetry we may assume that $a_1$, $a_2$, $a_3$, and $a_4$ have colors in $\{ 1, 2, 3 \}$.
	Next we proceed only with $H$ to obtain the contradiction.
	
	By symmetry, assume $b_4$ is the unique vertex in $H$ that (possibly) has color $4$.
	Without loss of generality, we assume $a_1$, $b_1$, and $a_2$ are colored by $x$, $y$, and $z$, respectively,
	where $\{x,y,z\} = \{1,2,3\}$.
	This forces $b_2$ to be colored with $x$, $a_3$ to be colored with $y$, and $b_3$ to be colored with $z$.
	Since $a_4$ is adjacent to vertices with colors $x$, $y$, and $z$, it must have color $4$, a contradiction.
\end{proof}

\begin{figure}
\[
\begin{array}{cc}
\vc{
	\includegraphics{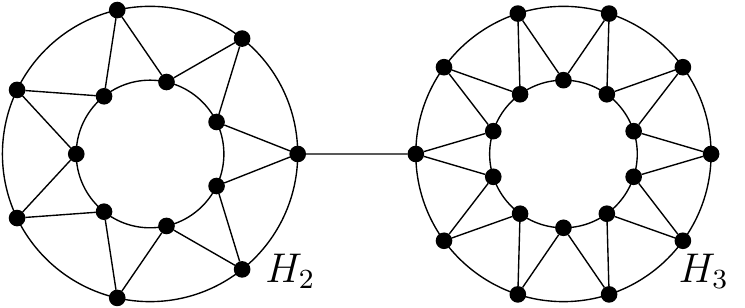}
}
&
\vc{
	\includegraphics{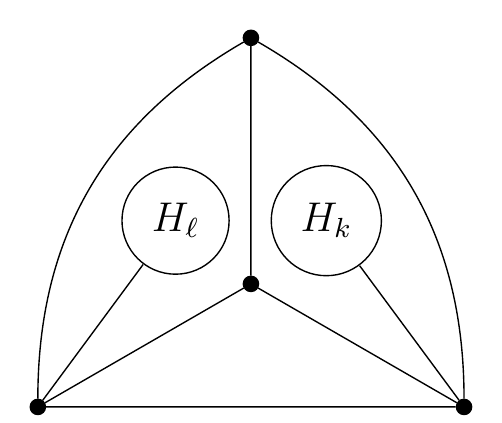}
}
\\[20pt]
 \text{(a)} & \text{(b)} 
\end{array}
\]
	\caption{More counterexamples to Conjecture~\ref{conj:plane4}.}
	\label{fig:plane42}
\end{figure}

The contradiction in Proposition~\ref{prop} is produced from the property of $H$ that every coloring of $H$ by colors $\{1,2,3,4\}$,
where every interior face has a unique-maximum color, has a vertex in the outer face colored by $4$.
We can generalize the counterexample in Figure~\ref{fig:plane41} by constructing an infinite family
of graphs $\mathcal{H} = \{ H_k \}_{k \geq 1}$ that can take the place of $H$.
We construct a graph $H_k$ on $6k + 2$ vertices by first embedding the cycle
$b_1 b_2 \cdots b_{3k+1}$ inside the cycle $a_1 a_2 \cdots a_{3k+1}$.
For $1 \leq i \leq 3k$, add edges $a_i b_i$ and $b_i a_{i+1}$, then add the edges $a_{3k+1} b_{3k+1}$ and $b_{3k+1} a_1$.
By this definition, the graph $H$ is equivalent to $H_1$.
See Figure~\ref{fig:plane42}(a) for an example of a generalization of the counterexample.

It is possible to construct more diverse counterexamples by embedding copies of members of $\mathcal{H}$ inside the faces
of any $4$-chromatic graph $G$ and adding an edge from each copy to some vertex on the face it belongs to.
It suffices to embed the graphs from $\mathcal{H}$ into a set of faces $K$ such that in every $4$-coloring of $G$,  there is at least one face in $K$ incident with a vertex of $G$ colored by 4.
An example of this with $G$ being $K_4$ is given in Figure~\ref{fig:plane42}(b).

We now introduce a variation of Conjecture~\ref{conj:plane4} with maximum degree and connectivity conditions added.

\begin{conjecture}
	\label{conj:plane4new}
	If $G$ is a connected plane graph with maximum degree $4$, then $\chiC (G) \leq 4$.
\end{conjecture}

Notice that we constructed a counterexample of maximum degree five.
Moreover, removing the edge $a_4 a_2'$ from the graph in Figure~\ref{fig:plane41} gives a disconnected graph
with maximum degree $4$ that does not have a FUM-coloring with colors in $\{1, 2, 3, 4\}$.
Recall that Andova et al.~\cite{1708.00094} showed that maximum degree $3$ suffices.

For a surface $\Sigma$, we define the facial unique-maximum chromatic number of $\Sigma$,
\[
\chiC(\Sigma) = \max_{ G\hookrightarrow \Sigma} \ \chiC(G),
\]
as the maximum of  $ \chiC(G)$ over all graphs $G$ embedded into $\Sigma$.
Our construction and the result of Wendland~\cite{Wen16} implies that 
$\chiC(S_0) = 5$, where $S_0$ is the sphere.
Our result motivates to study this invariant for graphs on other surfaces.
It would be interesting to have a similar characterization to Heawood number for other surfaces of higher genus.

\begin{problem}
Determine $\chiC(\Sigma)$ for surfaces $\Sigma$ of higher genus.
\end{problem}

\paragraph{Acknowledgment.}

The project has been supported by the bilateral cooperation between USA and Slovenia, project no. BI--US/17--18--013.
R. \v{S}krekovski was partially supported by the Slovenian Research Agency Program P1--0383.
B. Lidick\'y was partially supported by NSF grant DMS-1600390.

\bibliographystyle{abbrv}
\bibliography{MainBase}

\begin{thebibliography}{1}

\bibitem{1708.00094}
V.~Andova, B.~Lidick\'y, B.~Lu\v{z}ar, and R.~\v{S}krekovski.
\newblock On facial unique-maximum (edge-)coloring, 2017.
\newblock submitted.

\bibitem{AppHak76}
K.~Appel and W.~Haken.
\newblock {The solution of the four-color map problem}.
\newblock {\em Sci. Amer.}, 237:108--121, 1977.

\bibitem{CzaJen16}
J.~Czap and S.~Jendrol'.
\newblock {Facially-constrained colorings of plane graphs: A survey}.
\newblock {\em Discrete Math.}, 2016.
\newblock published online.

\bibitem{FabGor16}
I.~Fabrici and F.~G\"{o}ring.
\newblock {Unique-maximum coloring of plane graphs}.
\newblock {\em Discuss. Math. Graph Theory}, 36(1):95, 2016.

\bibitem{FabJenVrb15}
I.~Fabrici, S.~Jendrol', and M.~Vrbjarov\'{a}.
\newblock {Unique-maximum edge-colouring of plane graphs with respect to
  faces}.
\newblock {\em Discrete Appl. Math.}, 185:239--243, 2015.

\bibitem{Wen16}
A.~Wendland.
\newblock {Coloring of Plane Graphs with Unique Maximal Colors on Faces}.
\newblock {\em J. Graph Theory}, 83(4):359--371, 2016.

\end{thebibliography}

\end{document}